\documentclass[11pt]{amsart}
\usepackage{iftex}
\ifPDFTeX
  \usepackage[T1]{fontenc}
  \usepackage[utf8]{inputenc}
  \usepackage{lmodern}
\else
  \usepackage{fontspec}
\fi
\usepackage{amsmath,amssymb,amsthm,mathtools}
\ifPDFTeX\else
  \usepackage{unicode-math}
\fi
\usepackage[a4paper,margin=28mm]{geometry}
\usepackage{microtype}
\usepackage{xcolor}
\usepackage{cite}
\usepackage[colorlinks=true,linkcolor=blue!45!black,citecolor=blue!45!black,urlcolor=blue!45!black]{hyperref}
\numberwithin{equation}{section}
\newtheorem{theorem}{Theorem}[section]
\newtheorem{proposition}[theorem]{Proposition}
\newtheorem{lemma}[theorem]{Lemma}

\theoremstyle{remark}
\newcommand{\Rea}{\operatorname{Re}}
\newcommand{\Ima}{\operatorname{Im}}
\newcommand{\Res}{\operatorname*{Res}}
\newcommand{\eps}{\varepsilon}
\newcommand{\Z}{\mathbb Z}
\newcommand{\Q}{\mathbb Q}

\newcommand{\N}{\mathrm N}
\newcommand{\OO}{\mathcal O}
\newcommand{\dd}{\,\mathrm d}
\newcommand{\chip}[2]{\left(\frac{#1}{#2}\right)_4}
\newcommand{\zk}{\zeta_K^{(2)}}
\newcommand{\sumstar}{\sideset{}{^{*}}{\sum}}
\title[First moments of quartic and sextic Dirichlet $L$-functions]{Improved first moments of primitive quartic and sextic Dirichlet $L$-functions}
\author{Ziwei Hong}
\address{Renmin University of China, Beijing 100872, China}
\email{hongziwei@live.com}
\address{UNSW Sydney, Sydney NSW 2033, Australia}
\email{ziwei.hong@unsw.edu.au}
\date{}
\begin{document}
\begin{abstract}
We study the smoothed first moments of primitive quartic and sextic Dirichlet $L$-functions for which every nontrivial power of the character remains primitive. Assuming the generalized Riemann hypothesis for the Dirichlet and finite order Hecke $L$-functions occurring in the proof, we obtain asymptotic formulas with errors $O(Q^{\frac{49}{58}+\eps})$ and $O(Q^{\frac{23}{26}+\eps})$, respectively. We follow the double Dirichlet series approach of Gao and Zhao. The main ingredient is a joint estimate for the two sums arising from the removal of rational prime divisors, with the M\"obius coefficients retained. This gives a quotient of Dirichlet $L$-functions. Under GRH, we bound both the Dirichlet $L$-functions and their reciprocals. Combining the joint estimate with estimates for the generating series of Gauss sums gives larger regions of meromorphic continuation. These regions do not separate the possible secondary terms.
\end{abstract}
\maketitle

\section{Introduction}
\subsection{Background}
We consider the first moments of central values of primitive quartic and sextic Dirichlet $L$-functions. Gao and Zhao obtain asymptotic formulas for these families under the generalized Lindel\"of hypothesis, with errors $O(Q^{\frac{9}{10}+\eps})$ and $O(Q^{\frac{13}{14}+\eps})$, respectively; see [Theorem 1.1 of~\citen{GZ}]. Their argument uses double Dirichlet series and the analytic properties of generating series associated with Gauss sums. We follow this approach and obtain smaller errors under the generalized Riemann hypothesis.

Using multiple Dirichlet series, \v{C}ech~\cite{CechSecondary} separates a secondary term in the first moment of quadratic Dirichlet $L$-functions under the Riemann hypothesis and the generalized Lindel\"of hypothesis, and obtains continuation results and estimates for higher moments of the quadratic family in~\cite{CechMoments}. The arithmetic of quadratic characters and Gauss sums is different from that of the characters considered here, so these arguments do not apply directly.

Diaconu, Ion, Pa\c{s}ol and Popa~\cite{DIPP} develop a multiple Dirichlet series approach to the twisted first and second moments of Hecke $L$-functions associated with characters of order $r\geq3$. Their general setting assumes that the base field contains the $2r$th roots of unity. They obtain moment asymptotics for families indexed by squarefree ideals and by $r$th power free ideals, and combine the first and second mollified moments to prove that a positive proportion of the central values are nonzero. For $r=4$, their first moment error exponents are $\frac{11}{14}$ and $\frac{7}{11}$ in these two families. For $r=6$, the corresponding exponents are $\frac{5}{6}$ and $\frac{11}{17}$; see [Proposition 10.9.2 and Theorem 10.11.2 of~\citen{DIPP}]. Here squarefree describes the parameter ideal, not the order of the character. These results concern partial Hecke $L$-functions averaged over ideals by norm, whereas we consider Dirichlet characters over $\Q$ averaged by conductor. In particular, our errors are not asserted to improve the errors in~\cite{DIPP}.

Castillo, de Faveri and Dunn~\cite{CDF} study moments and nonvanishing for quartic Hecke $L$-functions over the Gaussian field. Their first moment result concerns a Hecke family, rather than the rational Dirichlet family considered here. That work does not give a sextic first moment asymptotic. We use estimates of David, Dunn, Hamieh and Lin~\cite{DDHL} for quartic Gauss sums, but no extrapolation of their quartic estimates is used for the sextic case.

\subsection{The families and the main results}
Let $j\in\{4,6\}$ and let $\chi_0$ denote the principal Dirichlet character. The notation $\sumstar$ in a sum with the condition $\chi^j=\chi_0$ denotes characters of exact order $j$ such that every character $\chi^k$, $1\leq k<j$, is primitive modulo $q$. For $j=4$, this is equivalent to requiring that $\chi$ and $\chi^2$ are primitive. For $j=6$, all five nontrivial powers are required to remain primitive; it is not enough that the global order of $\chi$ be six. We always impose $(q,j)=1$ and omit $q=1$. Such characters exist precisely when $q>1$ is squarefree and every prime divisor of $q$ is congruent to $1$ modulo $j$.

Let $\Phi\in C_c^\infty((0,\infty))$ and write
\begin{equation}\label{eq:mellin-transform}
 \widehat\Phi(s)=\int_0^\infty\Phi(x)x^s\frac{\dd x}{x}.
\end{equation}
Put $K_4=\Q(i)$, $K_6=\Q(\sqrt{-3})$, and
\begin{equation}\label{eq:error-exponents}
 \theta_4=\frac{49}{58},\qquad \theta_6=\frac{23}{26}.
\end{equation}
Our main result is as follows.

\begin{theorem}\label{thm:main}
Let $j\in\{4,6\}$. Assume the generalized Riemann hypothesis for Dirichlet $L$-functions and for finite order Hecke $L$-functions over $K_j$. For every $\Phi\in C_c^\infty((0,\infty))$ and every $\eps>0$, we have
\begin{equation}\label{eq:main}
 \sum_{\substack{q>1\\(q,j)=1}}\ \sumstar_{\substack{\chi\bmod q\\\chi^j=\chi_0}}
 L(\tfrac12,\chi)\Phi\left(\frac{q}{Q}\right)
 =C_j\widehat\Phi(1)Q+O_{\Phi,\eps}(Q^{\theta_j+\eps}),
\end{equation}
where $C_4>0$ and $C_6>0$ are given in \eqref{eq:constant} and \eqref{eq:sextic-constant}, respectively.
\end{theorem}

The gains in the exponents over [Theorem 1.1 of~\citen{GZ}] are
\begin{equation}\label{eq:exponent-comparison}
 \frac9{10}-\frac{49}{58}=\frac8{145},\qquad
 \frac{13}{14}-\frac{23}{26}=\frac4{91}.
\end{equation}
The hypotheses are stronger here: we use the generalized Riemann hypothesis to bound both $L$-functions and the factors $1/L$ occurring in the proof. It suffices to assume it for the characters and their powers used below. The constants $C_j$ are the same first main coefficients as in~\cite{GZ}; the improvement concerns the errors.

\subsection{The argument}
We follow the notation and the double Dirichlet series framework of Gao and Zhao~\cite{GZ}. The main difference is that we retain the M\"obius function in the sum imposing the restriction on rational prime divisors, rather than estimating this sum directly by absolute values. After interchanging the sums, we estimate the sums over $m$ and $d$ together and obtain a quotient of Dirichlet $L$-functions. The M\"obius coefficients produce the $L$-function in the denominator. Under GRH, both Mellin variables can be moved to the right of $\frac12$. The denominator involves the square of a quartic character in the quartic case and the fourth power of a sextic character in the sextic case.

For $j=4$, we combine this estimate with [(4.34) of~\citen{DDHL}], retaining the factor $\zeta_{K_4}^{(2)}(4u-3)$. Interpolation and Mellin inversion give a region whose relevant boundary is $58\Rea s+40\Rea w=69$. For $j=6$, [Lemma 2.9 of~\citen{GZ}] gives the required estimate for the coprime sextic Gauss series in $\Rea u>1$. The rational divisor occurs in the frequency as $md^4$, rather than $md^2$, and the corresponding boundary is $26\Rea s+20\Rea w=33$.

We hoped to use this approach to separate the possible secondary terms. The possible second poles lie on $s+w=\frac54$ and $s+w=\frac76$, which meet $w=\frac12$ at $s=\frac34$ and $s=\frac23$, respectively. Both points remain to the left of the contours obtained here. Further analytic input is therefore still needed to separate these terms.

Section~\ref{sec:preliminaries} records the quartic notation, the estimates for Gauss sum generating series, and the results on functions of several complex variables used in both cases. Section~\ref{sec:proof} proves the quartic case of Theorem~\ref{thm:main}. Section~\ref{sec:sextic} establishes the corresponding sextic estimates and completes the proof. Throughout the paper, $\eps$ denotes a positive number which may change from one occurrence to the next.

\section{Preliminaries}\label{sec:preliminaries}
\subsection{Quartic residue symbols and characters}
In this section and Section~\ref{sec:proof}, let $K=K_4=\Q(i)$, $\OO_K=\Z[i]$ and $\lambda=1+i$. An odd Gaussian integer is called primary if it is congruent to $1$ modulo $\lambda^3$. Every nonzero ideal of $\Z[i]$ relatively prime to $2$ has a unique primary generator. We write $|z|$ for the modulus of a complex number $z$. For $z\in\OO_K$, its norm is $\N(z)=z\bar z=|z|^2$; for a nonzero ideal $\mathfrak a\subseteq\OO_K$, we write $\N(\mathfrak a)=\#(\OO_K/\mathfrak a)$. The quartic residue symbol is extended by zero when its arguments are not coprime.

We write
\begin{equation}\label{eq:characters}
 \chi_4^{(m)}(c)=\chip{m}{c},\qquad \chi_n(x)=\chip{x}{n}.
\end{equation}
The first symbol is regarded as a Hecke character on odd ideals, using their primary generators. A character in the quartic family of Theorem~\ref{thm:main} can be written uniquely as
\begin{equation}\label{eq:param}
 \chi_n(m)=\chip{m}{n},\qquad q=\N(n),
\end{equation}
where $n$ is primary, squarefree in $\Z[i]$, has no rational prime divisor, and is not $1$. This is the case $j=4$ of [Lemma 2.2 of~\citen{GZ}]. In particular, a split prime $p=\varpi\bar\varpi$ contributes at most one of $\varpi$ and $\bar\varpi$ to $n$.

For coprime primary $c,n$, quartic reciprocity gives
\begin{equation}\label{eq:reciprocity}
 \chip{c}{n}=\chip{n}{c}
 (-1)^{\frac{\N(c)-1}{4}\frac{\N(n)-1}{4}}.
\end{equation}
For coprime rational integers $m,d$ with $d$ odd, $\chip{m}{d}=1$; see [Lemma 2.3 of~\citen{GZ}]. We use the primary associate of a rational Gaussian integer whenever it occurs as a modulus.

\subsection{Gauss sums}
We use
\begin{equation}\label{eq:gauss}
 e_K(z)=\exp(2\pi i\Ima z),\qquad
 g_K(r,\chi_n)=\sum_{x\bmod n}\chip{x}{n}e_K\left(\frac{rx}{n}\right).
\end{equation}
We abbreviate $g_K(1,\chi_n)$ to $g_K(\chi_n)$.
For a primitive Dirichlet character $\chi$ of conductor $q$, write
\begin{equation}\label{eq:tau-parity}
 \tau(\chi)=\sum_{a\bmod q}\chi(a)\exp\left(\frac{2\pi ia}{q}\right),
 \qquad \chi(-1)=(-1)^\delta,\quad \delta\in\{0,1\}.
\end{equation}
With the convention \eqref{eq:gauss}, the relation needed here is
\begin{equation}\label{eq:gauss-normalization}
 g_K(\chi_n)=i^{-\delta}\tau(\chi_n).
\end{equation}
To specify the normalization, the change of additive character gives
\begin{equation}\label{eq:tau-change}
\tau(\chi_n)=\overline{\chip{2i}{n}}\chip{\bar n}{n}g_K(\chi_n),
\end{equation}
whereas the supplementary law gives
\begin{equation}\label{eq:supplementary}
 \chip{\bar n}{n}=i^\delta\chip{2i}{n}.
\end{equation}
Combining \eqref{eq:tau-change} and \eqref{eq:supplementary}, the factors involving $\chip{2i}{n}$ therefore cancel. We use this identity with the additive character defined in \eqref{eq:gauss}.
The trace convention in \cite{DDHL} is related to ours by
\begin{equation}\label{eq:trace}
 g_K(r,\chi_n)=g_4\left(\frac{r}{2i},n\right).
\end{equation}
We keep this change of additive character throughout.

For the characters in \eqref{eq:param}, define
\begin{equation}\label{eq:beta-parity}
 \beta_0=1,\qquad \beta_1=1+\lambda^3=-1+2i\pmod4,
 \qquad n\equiv\beta_\delta\pmod4
 \ \Longleftrightarrow\ \chi_n(-1)=(-1)^\delta.
\end{equation}
Thus $\delta=0$ corresponds to an even character and $\delta=1$ to an odd character. The parameter $\beta$ below belongs to $\{\beta_0,\beta_1\}$.

\subsection{Estimates for \texorpdfstring{$L$}{L}-functions}
The functional equation and Stirling estimates used below are standard; see \cite{IK}. Under the hypothesis of Theorem~\ref{thm:main}, we use Lindel\"of bounds for $L(z,\chi)$ and $\frac{1}{L(z,\chi)}$, together with the analogous Hecke $L$-function bounds, in each fixed half plane $\Rea z\geq\frac{1}{2}+\eta$, where $\eta>0$. These estimates are uniform in the analytic conductor. A pole of a principal $L$-function is removed before applying the corresponding bound. For a principal character $\chi$, the function $\frac{1}{L(z,\chi)}$ has a zero at $z=1$.

For a primitive Dirichlet character $\chi$ modulo $q$ with $\chi(-1)=(-1)^\delta$, the functional equation is
\begin{equation}\label{eq:dirichlet-FE}
 L(w,\chi)=\frac{i^{-\delta}\tau(\chi)}{q^w}
 \gamma_\delta(w)L(1-w,\bar\chi),\qquad
 \gamma_\delta(w)=\pi^{w-\frac{1}{2}}
 \frac{\Gamma\left(\frac{1-w+\delta}{2}\right)}{\Gamma\left(\frac{w+\delta}{2}\right)}.
\end{equation}
On fixed real strips, Stirling's formula and the hypothesis give
\begin{equation}\label{eq:conductor-bound}
 L(w,\chi)\ll_{\eps}q^{\max(\frac{1}{2}-\Rea w,0)+\eps}
 (1+|w|)^B
\end{equation}
for nonprincipal primitive $\chi$, with $B$ depending on the strip. We also use the analogous bounds for Hecke $L$-functions over $K$. Euler factors omitted from imprimitive $L$-functions cost an arbitrarily small power of the modulus on the half planes used below.

\subsection{The generating series of quartic Gauss sums}
For $\beta\in\{\beta_0,\beta_1\}$ and a positive integer $r$, set
\begin{equation}\label{eq:H}
 H_\beta(u,r)=\sum_{\substack{c\equiv\beta\ (4)\\(c,r)=1}}
                    \frac{g_K(r,\chi_c)}{\N(c)^u}.
\end{equation}
The summation is over nonzero Gaussian integers in the indicated primary class. If $(c,r)=1$, then $g_K(r,\chi_c)=\overline{\chip{r}{c}}g_K(\chi_c)$ and the Gauss sum vanishes unless $c$ is squarefree. Thus no extra squarefree condition is required in \eqref{eq:H}.

Write
\begin{equation}\label{eq:zeta-omitted}
 \zk(z)=\zeta_K(z)(1-2^{-z}).
\end{equation}
We retain the factor in \eqref{eq:zeta-omitted} in the following estimate. The bound for the completed series is available in a larger region than the bound obtained by dividing by the zeta factor.

\begin{proposition}\label{prop:coprime}
The function $(u-\frac54)\zk(4u-3)H_\beta(u,r)$ is holomorphic for $\Rea u>\frac{3}{4}$. Given $\eta>0$ and $\eps>0$, uniformly for $\frac{3}{4}+\eta\leq\sigma=\Rea u\leq\frac{3}{2}$,
\begin{equation}\label{eq:coprime}
 (u-\frac54)\zk(4u-3)H_\beta(u,r)\ll_{\eta,\eps}
 r^{\frac{3}{2}-\sigma+\eps}(1+|u|)^B,
\end{equation}
where $B$ is independent of $r$.
\end{proposition}
\begin{proof}
We give the local argument because uniformity in the coprimality condition is needed later. Let $\psi_i(u,\nu)$ be the unrestricted Gauss sum series in the finite system of [Section 4 of~\citen{DDHL}], in their trace convention, and put
\begin{equation}\label{eq:completed-gauss}
 Z_i(u,\nu)=\zk(4u-3)\psi_i(u,\nu).
\end{equation}
Equations (4.31)--(4.34) of that paper give the twist exponent $\frac{\frac{3}{2}-\sigma}{2}$ for $Z_i$. Apply the same Phragm\'en--Lindel\"of argument after multiplying by $(u-\frac{3}{4})(u-\frac{5}{4})$. Multiplication by this polynomial changes the vertical bounds only by a polynomial in $u$. Dividing by $u-\frac{3}{4}$ on $\sigma\geq\frac{3}{4}+\eta$ gives
\begin{equation}\label{eq:full-bound}
 (u-\frac54)Z_i(u,\nu)
 \ll_{\eta,\eps,\operatorname{ord}_\lambda\nu}
 \N(\nu)^{\frac{\frac{3}{2}-\sigma}{2}+\eps}(1+|u|)^B.
\end{equation}
This also covers a neighborhood of $u=\frac{5}{4}$, rather than using an estimate only outside a disk around that point.

Let $\varpi$ be an odd primary Gaussian prime and put $q=\N(\varpi)$. Fix $\beta\in\{\beta_0,\beta_1\}$ and a nonzero $t\in\lambda^{-2}\OO_K$ with $(t,\varpi)=1$. Let $d$ be the product of the odd primary Gaussian primes for which coprimality has already been imposed, with $(d,\varpi)=1$ and $d=1$ at the first step. Initially for $\Rea u>2$, define
\begin{equation}\label{eq:local-series-definitions}
 \begin{aligned}
 F_0&=\sum_{\substack{c\equiv\beta\ (4)\\(c,d)=1}}
       \frac{g_4(t,c)}{\N(c)^u},&
 F_2&=\sum_{\substack{c\equiv\varpi\beta\ (4)\\(c,d)=1}}
       \frac{g_4(t\varpi^2,c)}{\N(c)^u},\\
 X&=\sum_{\substack{c\equiv\beta\ (4)\\(c,d\varpi)=1}}
       \frac{g_4(t,c)}{\N(c)^u},&
 Y&=\sum_{\substack{c\equiv\varpi\beta\ (4)\\(c,d\varpi)=1}}
       \frac{g_4(t\varpi^2,c)}{\N(c)^u}.
 \end{aligned}
\end{equation}
All four sums are over nonzero $c\in\OO_K$, and $g_4$ is the Gauss sum in the trace convention of~\cite{DDHL}. The dependence on $u,t,\beta,d,\varpi$ is suppressed in these temporary notations. Thus $F_0,F_2$ include terms with $\varpi\mid c$, whereas $X,Y$ satisfy the additional condition $(c,\varpi)=1$; the condition $(c,d)=1$ is retained in all four sums. These definitions do not include the common factor $\zk(4u-3)$. We now give the relations between $F_0,F_2$ and $X,Y$. Separating the terms according to the power of $\varpi$ dividing $c$ and applying the local Gauss sum identities gives
\begin{equation}\label{eq:system}
 F_0=X+aY,\qquad F_2=Y+bX,\qquad
 ab=q^{3-4u},\quad |a|=q^{\frac{1}{2}-\sigma},\quad |b|=q^{\frac{5}{2}-3\sigma}.
\end{equation}
With the congruence classes specified in \eqref{eq:local-series-definitions}, one can take
\begin{equation}\label{eq:local-coefficients}
 \begin{aligned}
 a&=q^{-u}(-1)^{C(\varpi,\varpi\beta)}g_4(t,\varpi),\\
 b&=q^{2-3u}\chip{-1}{\varpi}(-1)^{C(\varpi,\beta)}\overline{g_4(t,\varpi)},\\
 C(x,y)&=\frac{\N(x)-1}{4}\frac{\N(y)-1}{4}.
 \end{aligned}
\end{equation}
The same relations hold after multiplying all four series by $\zk(4u-3)$. Solving for $X,Y$, we obtain
\begin{equation}\label{eq:local-inversion}
 X=\frac{F_0-aF_2}{1-q^{3-4u}},\qquad
 Y=\frac{F_2-bF_0}{1-q^{3-4u}}.
\end{equation}
For twists $t\varpi$ and $t\varpi^3$, the relations between the series with and without the condition $(c,\varpi)=1$ give the respective factors
\begin{equation}\label{eq:odd-states}
 \frac1{1+q^{1-2u}g_2(t,\varpi)},\qquad
 \frac1{1-q^{3-4u}},\qquad |g_2(t,\varpi)|=q^{\frac{1}{2}}.
\end{equation}
Here $g_2$ is the quadratic Gauss sum in the trace convention of \cite{DDHL}. These identities follow by separating the allowed valuations in the local formulas of [Section 3 of~\citen{DDHL}]; all changes of primary class are retained. Every denominator in \eqref{eq:system}--\eqref{eq:odd-states} is nonzero for $\sigma>\frac{3}{4}$.

Put $a_0=\frac{\frac{3}{2}-\sigma}{2}\geq0$. Suppose the original twist has valuation at $\varpi$ $k=4h+j$, $0\leq j\leq3$. In the series with the coprimality condition we may first remove the fourth power, but must retain the condition $(c,\varpi)=1$ if $k>0$. Compare every resulting term with the original bound $q^{ka_0}$ in \eqref{eq:full-bound}.

If $j=0$ and $h\geq1$, the two terms in \eqref{eq:local-inversion} for $X$ have relative sizes
\begin{equation}\label{eq:local-zero}
 q^{-4ha_0},\qquad q^{\frac{1}{2}-\sigma+(2-4h)a_0}.
\end{equation}
The first expression in \eqref{eq:local-zero} is at most $1$, and the second is at most $q^{-1}$, since for $h=1$ its exponent is exactly $-1$. If $j=2$, the two terms in the formula for $Y$ have relative sizes at most $1$ and
\begin{equation}\label{eq:local-two}
 q^{\frac{5}{2}-3\sigma-(4h+2)a_0}\leq q^{1-2\sigma}.
\end{equation}
For $j=1,3$, removing the fourth power only decreases the bound. Thus the factor introduced by imposing $(c,\varpi)=1$ is bounded, respectively, by
\begin{equation}\label{eq:costs}
 \frac{1+q^{-1}}{1-q^{3-4\sigma}},\quad
 \frac1{1-q^{\frac{3}{2}-2\sigma}},\quad
 \frac{1+q^{1-2\sigma}}{1-q^{3-4\sigma}},\quad
 \frac1{1-q^{3-4\sigma}}.
\end{equation}
Each expression is $1+O_\eta(q^{-c_\eta})$ for some $c_\eta>0$. Their product over the prime divisors of a twist is $O_{\eta,\eps}(\N(r)^\eps)$: large primes are absorbed into $q^\eps$, and there are only finitely many smaller primes. At each prime, the displayed bound includes both terms in the local relation. Iteration therefore bounds the full sum obtained by imposing $(c,\varpi)=1$ for every prime $\varpi$ dividing the twist.

For rational $r$, both primes above a split rational prime must be included. Since $\N(r)=r^2$, the exponent in \eqref{eq:full-bound} becomes $r^{\frac{3}{2}-\sigma}$. The intermediate twists need not be rational; \eqref{eq:full-bound} applies to their Gaussian norms. Finally $H_\beta(u,16r)=H_\beta(u,r)$, so the valuation at $\lambda$ reduces to finitely many cases. This proves \eqref{eq:coprime}. The same finite relations, together with the holomorphic continuation of $(u-\frac34)(u-\frac54)Z_i$, prove the asserted holomorphy.
\end{proof}

\subsection{Multiple complex variables}
For an open set $U\subset\mathbb R^n$, define the associated tube domain by
\begin{equation}\label{eq:tube-domain}
 T(U)=\{z\in\mathbb C^n:\Rea z\in U\}=U+i\mathbb R^n.
\end{equation}
We write $\widehat T$ for the convex hull of a tube domain $T$. We first recall Bochner's tube theorem \cite{Bochner}, in the form stated in [Theorem 2.12 of~\citen{GZ}].

\begin{theorem}[Bochner]\label{thm:tube}
Suppose that $U\subset\mathbb R^n$ is open and connected. Every function holomorphic on $T(U)$ has a unique holomorphic continuation to $\widehat{T(U)}$.
\end{theorem}

The next result is [Proposition 2.13 of~\citen{GZ}]. It allows bounds on a tube domain to be retained after continuation to its convex hull.

\begin{samepage}
\begin{theorem}\label{thm:tube-bound}
Let $T\subset\mathbb C^n$ be a connected tube domain. Suppose that $g$ and $h$ are holomorphic on $T$ and that $h$ has no zeros there. Denote their holomorphic continuations to $\widehat T$ by $\widetilde g$ and $\widetilde h$. If
\begin{equation}\label{eq:tube-bound-source}
 |g(z)|\leq|h(z)|\qquad(z\in T),
\end{equation}
then
\begin{equation}\label{eq:tube-bound-hull}
 |\widetilde g(z)|\leq|\widetilde h(z)|\qquad(z\in\widehat T).
\end{equation}
\end{theorem}
\end{samepage}

In each application of Theorems~\ref{thm:tube} and \ref{thm:tube-bound}, we first multiply by the factors needed to obtain a holomorphic function.

\section{The quartic first moment}\label{sec:proof}
\subsection{Initial treatment}
For $\Rea s$ and $\Rea w$ sufficiently large, define
\begin{equation}\label{eq:A}
 A_4(s,w)=\sum_{\substack{q>1\\(q,2)=1}}\ \sumstar_{\substack{\chi\bmod q\\\chi^4=\chi_0}}\frac{L(w,\chi)}{q^s}
 =\sideset{}{'}{\sum}_{n\ {\rm primary}}
   \frac{L(w,\chi_n)}{\N(n)^s}.
\end{equation}
Here $\sum'$ restricts the sum to squarefree Gaussian integers $n\ne1$ with no rational prime divisor. This equality follows from \eqref{eq:param}. Mellin inversion gives
\begin{equation}\label{eq:moment-integral}
 \sum_{\substack{q>1\\(q,2)=1}}\ \sumstar_{\substack{\chi\bmod q\\\chi^4=\chi_0}} L(w,\chi)\Phi\left(\frac{q}{Q}\right)
 =\frac1{2\pi i}\int_{(c)}A_4(s,w)Q^s\widehat\Phi(s)\dd s,
\end{equation}
where $c$ is sufficiently large. We study the analytic properties of $A_4$ before shifting this contour at $w=\frac{1}{2}$.

\subsection{First region of convergence of \texorpdfstring{$A_4(s,w)$}{A4(s,w)}}
We begin in $\Rea s>1$ and $\Rea w>1$, where all the sums below converge absolutely. Expanding $L(w,\chi_n)$ and interchanging the sums in \eqref{eq:A}, we obtain
\begin{equation}\label{eq:first-expansion}
 A_4(s,w)=\sum_{m\geq1}\frac1{m^w}
 \sideset{}{'}{\sum}_{n\ {\rm primary}}
 \frac{\chi_4^{(m)}(n)}{\N(n)^s}.
\end{equation}
The prime on the sum has the same meaning as in \eqref{eq:A}, in particular $n=1$ is omitted. We first add this term back. Its contribution to the outer sum is $\zeta(w)$.

For a squarefree odd Gaussian integer $n$, the condition that $n$ has no rational prime divisor is detected by
\begin{equation}\label{eq:mobius-indicator}
 \mathbf1_{\substack{n\text{ has no rational}\\\text{prime divisor}}}
 =\sum_{\substack{d\mid n\\d\in\Z_{>0}}}\mu(d).
\end{equation}
Write $n=dc$, using the primary associate of $d$. Since $d$ is odd and squarefree as a rational integer, it is squarefree in $\OO_K$. Thus $dc$ is squarefree if and only if $c$ is squarefree and $(c,d)=1$. Moreover, $\N(d)=d^2$ and
\begin{equation}\label{eq:rational-symbol}
 \chi_4^{(m)}(d)=\chip{m}{d}
 =\begin{cases}1,&(m,d)=1,\\0,&(m,d)>1.\end{cases}
\end{equation}
Consequently, M\"obius inversion gives
\begin{equation}\label{eq:first-mobius}
 \begin{aligned}
 A_4(s,w)+\zeta(w)
 &=\sum_{m\geq1}\frac1{m^w}
 \sum_{\substack{d\geq1\\(d,2m)=1}}\frac{\mu(d)}{d^{2s}}
 \sum_{\substack{c\ {\rm primary}\text{ and squarefree}\\(c,2d)=1}}
 \frac{\chi_4^{(m)}(c)}{\N(c)^s}.
 \end{aligned}
\end{equation}

For any Hecke character $\psi$, let $L_K^{(2)}(z,\psi)$ denote its Euler product with the factor above $2$ omitted. For $\psi=\chi_4^{(m)}$, we retain its extension by zero at primes dividing $m$, including when it is imprimitive. Since each Gaussian prime can occur in a squarefree $c$ only to exponent zero or one, the innermost sum in \eqref{eq:first-mobius} has Euler product
\begin{align}
 &\sum_{\substack{c\ {\rm primary}\text{ and squarefree}\\(c,2d)=1}}
 \frac{\chi_4^{(m)}(c)}{\N(c)^s}
 =\prod_{\varpi\nmid2d}
 \left(1+\frac{\chi_4^{(m)}(\varpi)}{\N(\varpi)^s}\right)\notag\\
 &\qquad=
 \frac{L_K^{(2)}(s,\chi_4^{(m)})}
 {L_K^{(2)}(2s,(\chi_4^{(m)})^2)}
 \prod_{\varpi\mid d}
 \left(1+\frac{\chi_4^{(m)}(\varpi)}{\N(\varpi)^s}\right)^{-1}.
 \label{eq:squarefree-euler}
\end{align}
Indeed, the identity $1+x=\frac{1-x^2}{1-x}$ gives, at each Gaussian prime,
\begin{equation}\label{eq:squarefree-local}
 1+\frac{\chi_4^{(m)}(\varpi)}{\N(\varpi)^s}
 =\frac{1-(\chi_4^{(m)}(\varpi))^2\N(\varpi)^{-2s}}
 {1-\chi_4^{(m)}(\varpi)\N(\varpi)^{-s}}.
\end{equation}
The identity \eqref{eq:squarefree-local} explains both the argument $2s$ and the square of the character in the denominator of \eqref{eq:squarefree-euler}.

Substituting \eqref{eq:squarefree-euler} into \eqref{eq:first-mobius}, we obtain
\begin{equation}\label{eq:first-substitution}
 A_4(s,w)+\zeta(w)
 =\sum_{m\geq1}\frac1{m^w}
 \frac{L_K^{(2)}(s,\chi_4^{(m)})}
 {L_K^{(2)}(2s,(\chi_4^{(m)})^2)}
 \sum_{\substack{d\geq1\\(d,2m)=1}}\frac{\mu(d)}{d^{2s}}
 \prod_{\varpi\mid d}
 \left(1+\frac{\chi_4^{(m)}(\varpi)}{\N(\varpi)^s}\right)^{-1}.
\end{equation}

We now continue the representation in \eqref{eq:first-substitution}. Fix $\eta>0$ and write $\sigma=\Rea s\geq\frac{1}{2}+\eta$. The denominator in \eqref{eq:squarefree-euler} has an absolutely convergent Euler product, since $2\sigma>1$. The function $\frac{1}{L_K^{(2)}(2s,(\chi_4^{(m)})^2)}$ is holomorphic and bounded uniformly in $m$ in this half plane. Also,
\begin{equation}\label{eq:finite-euler-bound}
 \left|\prod_{\varpi\mid d}
 \left(1+\frac{\chi_4^{(m)}(\varpi)}{\N(\varpi)^s}\right)^{-1}\right|
 \leq\prod_{\varpi\mid d}(1-\N(\varpi)^{-\sigma})^{-1}
 \ll_{\eta,\eps}d^\eps.
\end{equation}
By \eqref{eq:finite-euler-bound}, choosing $0<\eps<2\eta$, we see that
\begin{equation}\label{eq:first-correction}
 \sum_{\substack{d\geq1\\(d,2m)=1}}
 \frac{\mu(d)}{d^{2s}}
 \prod_{\varpi\mid d}
 \left(1+\frac{\chi_4^{(m)}(\varpi)}{\N(\varpi)^s}\right)^{-1}
\end{equation}
converges normally for $\Rea s>\frac{1}{2}$ and is bounded uniformly in $m$ on fixed strips contained in $\Rea s>\frac{1}{2}$.

The conductor of $\chi_4^{(m)}$ is bounded by a fixed power of $m$. Under the stated hypothesis, the Hecke $L$-function bounds therefore give, on each fixed strip $\frac{1}{2}+\eta\leq\Rea s\leq A$,
\begin{equation}\label{eq:hecke-bound}
 (s-1)L_K^{(2)}(s,\chi_4^{(m)})
 \ll_{\eta,A,\eps}m^\eps(1+|s|)^B.
\end{equation}
The factor $s-1$ removes the possible pole when $\chi_4^{(m)}$ is principal. Combining \eqref{eq:first-substitution}, \eqref{eq:first-correction} and \eqref{eq:hecke-bound}, we obtain
\begin{equation}\label{eq:first-bound}
 |(s-1)(A_4(s,w)+\zeta(w))|
 \ll_{\eta,A,\eps}(1+|s|)^B
 \sum_{m\geq1}m^{-\Rea w+\eps}.
\end{equation}
For every $\Rea w>1$, choose $\eps<\Rea w-1$. The resulting series over $m$ converges normally on compact subsets, and $\zeta(w)$ is holomorphic there. We conclude that $(s-1)A_4(s,w)$ is holomorphic in
\begin{equation}\label{eq:initial-region}
 S_1=\{(s,w):\Rea s>\frac12,\ \Rea w>1\}.
\end{equation}
This argument follows the expansion in [Section 3.2 of~\citen{GZ}], with the square of the character retained in the squarefree Euler denominator. The convergence in $S_1$ concerns the series over $m$ after inserting \eqref{eq:squarefree-euler}; it supplies continuation of the original series from $\Rea s>1$, $\Rea w>1$.

\subsection{Second region of continuation of \texorpdfstring{$A_4(s,w)$}{A4(s,w)}}
We next apply \eqref{eq:conductor-bound} directly to the first representation of $A_4(s,w)$ in \eqref{eq:A}. Let $\omega(q)$ denote the number of distinct prime divisors of $q$. At every prime divisor of an admissible conductor there are two choices for the local quartic character. Thus, for every $\eps>0$,
\begin{equation}\label{eq:character-count}
 \sumstar_{\substack{\chi\bmod q\\\chi^4=\chi_0}}1
 \leq2^{\omega(q)}\ll_\eps q^\eps.
\end{equation}
Combining \eqref{eq:conductor-bound} and \eqref{eq:character-count}, and absorbing the two small powers of $q$ into $q^\eps$, we obtain, on each fixed real strip in $w$,
\begin{equation}\label{eq:conductor-absolute-bound}
 \sum_{\substack{q>1\\(q,2)=1}}
 \ \sumstar_{\substack{\chi\bmod q\\\chi^4=\chi_0}}
 \frac{|L(w,\chi)|}{q^{\Rea s}}
 \ll_\eps(1+|w|)^B
 \sum_{q\geq1}
 q^{-\Rea s+\max\{\frac{1}{2}-\Rea w,0\}+\eps}.
\end{equation}
The sum on the right of \eqref{eq:conductor-absolute-bound} converges after choosing $\eps$ sufficiently small whenever
\begin{equation}\label{eq:conductor-condition}
 \Rea s>1+\max\{\frac12-\Rea w,0\}
 =\max\{1,\frac32-\Rea w\}.
\end{equation}
For $\Rea w\geq\frac{1}{2}$, this is $\Rea s>1$. For $\Rea w<\frac{1}{2}$, it is $\Rea(s+w)>\frac{3}{2}$. Consequently, the first series in \eqref{eq:A} converges absolutely and normally on compact subsets of
\begin{equation}\label{eq:conductor-region}
 S_2=\{(s,w):\Rea s>1,\ \Rea(s+w)>\frac32\}.
\end{equation}
In particular, $A_4(s,w)$ is holomorphic in \eqref{eq:conductor-region}. This region connects the first representation to the functional equation used below.

\subsubsection{The functional equation and the M\"obius sum}
Let $A_{4,\delta}$ denote the part of \eqref{eq:A} with $\chi(-1)=(-1)^\delta$. Put
\begin{equation}\label{eq:C}
 C_{4,\delta}(u,v)=\sum_{\substack{q>1\\(q,2)=1}}
 \ \sumstar_{\substack{\chi\bmod q\\\chi^4=\chi_0\\\chi(-1)=(-1)^\delta}}
 \frac{i^{-\delta}\tau(\chi)L(v,\bar\chi)}{q^u}.
\end{equation}
Here the subscript $\delta\in\{0,1\}$ records parity. Thus $C_{4,\delta}$ is the parity decomposition of the reflected series; its index differs from the sign $a=\pm1$ used for $C_{4,a}$ in [Section 3.3 of~\citen{GZ}]. By \eqref{eq:dirichlet-FE},
\begin{equation}\label{eq:FE}
 A_{4,\delta}(s,w)=\gamma_\delta(w)C_{4,\delta}(s+w,1-w),
 \qquad A_4=A_{4,0}+A_{4,1}.
\end{equation}
\begin{lemma}\label{lem:mobius}
In a right hand domain of absolute convergence,
\begin{equation}\label{eq:mobius}
 C_{4,\delta}(u,v)+\mathbf1_{\delta=0}\zeta(v)
 =\sum_{m\geq1}\frac1{m^v}
   \sum_{\substack{d\geq1\\(d,2m)=1}}
   \frac{\mu(d)d}{d^{2u}}H_{\beta_\delta}(u,md^2).
\end{equation}
\end{lemma}
\begin{proof}
Use the Gaussian parametrization and insert $\sum_{d\mid n,\ d\in\Z_{>0},\ 2\nmid d}\mu(d)$ to exclude rational prime divisors. For odd squarefree rational $d$, use its primary associate. Quartic reciprocity, together with the fact that the quartic symbol between coprime rational integers is $1$, gives
\begin{equation}\label{eq:gauss-rational}
 g_K(\chi_d)=d,\qquad
 g_K(m,\chi_{dc})=d\,g_K(md^2,\chi_c)\quad ((m,d)=(c,d)=1).
\end{equation}
Using \eqref{eq:gauss-rational}, we obtain \eqref{eq:mobius}. The Gaussian sum includes $n=1$ in the even class; its contribution is precisely $\zeta(v)$, which explains the term on the left. In particular no additional quadratic character of $d$ occurs in this normalization.
\end{proof}

\subsubsection{Key step}
The next identity supplies cancellation in the estimate on the line $\Rea u=\frac{3}{2}+\eta$, where $\eta>0$. For a possibly imprimitive Dirichlet character $\chi$, extended by zero at the primes of its modulus, let $L^{(2)}$ denote the Euler product with the factor at $2$ omitted.

\begin{lemma}\label{lem:euler}
For $\Rea z,\Rea t>1$, we have
\begin{equation}\label{eq:euler}
 \sum_{\substack{m,d\geq1\\(m,d)=1,\ 2\nmid d}}
 \frac{\mu(d)\chi(m)\chi(d)^2}{m^z d^t}
 =\frac{L(z,\chi)}{L^{(2)}(t,\chi^2)}
 \prod_{p\ne2}\left(1+
 \frac{\chi(p)^3p^{-z-t}}{1-\chi(p)^2p^{-t}}\right).
\end{equation}
The last product converges normally when $\Rea t>0$ and $\Rea(z+t)>1$.
\end{lemma}
\begin{proof}
At an odd prime write $a=\chi(p)p^{-z}$ and $b=\chi(p)^2p^{-t}$. The local sum is $(1-a)^{-1}-b=\frac{1-b+ab}{1-a}$, which is the local factor on the right. At $2$, only the $m$-sum remains. On compact subsets of the stated region the correction differs from $1$ by $O(p^{-\Rea(z+t)})$, uniformly in $\chi$.
\end{proof}

Let $U(x)$ and $V(x)$ be fixed smooth functions compactly supported on the set of positive real numbers $\mathbb R_+$. We write $\widehat U$ and $\widehat V$ for their Mellin transforms, with the same convention as in \eqref{eq:mellin-transform}. Set
\begin{equation}\label{eq:weight}
 \begin{aligned}
 W_u(x)&=(x\partial_x+1)(x\partial_x+2-u)U(x),\\
 \widehat W_u(z)&=\int_0^\infty W_u(x)x^z\frac{\dd x}{x}
 =(1-z)(2-u-z)\widehat U(z).
 \end{aligned}
\end{equation}
The Mellin identity follows by two integrations by parts; compact support makes the boundary terms vanish. The first factor removes the principal character pole encountered in the Mellin argument. The second keeps the first polar plane explicit on reconstruction. The weights $W_u\left(\frac{m}{M}\right)$ and $V\left(\frac{d}{D}\right)$ restrict the sums to $m\asymp M$ and $d\asymp D$, respectively. Define
\begin{equation}\label{eq:weighted-sum}
 T_\beta(u;M,D)=(u-\frac54)\zk(4u-3)
 \sum_{\substack{m,d\geq1\\(d,2m)=1}}
 \mu(d)d\,H_\beta(u,md^2)W_u\left(\frac{m}{M}\right)V\left(\frac{d}{D}\right).
\end{equation}
The sums in \eqref{eq:weighted-sum} are finite for fixed $M,D$. By Proposition~\ref{prop:coprime}, $T_\beta(u;M,D)$ is holomorphic for $\Rea u>\frac{3}{4}$.

\begin{proposition}\label{prop:joint-bound}
Assume GRH as in Theorem~\ref{thm:main}. With $U,V$ as above, let $W_u$ and $T_\beta$ be given by \eqref{eq:weight} and \eqref{eq:weighted-sum}. Fix $0<\eta<\frac{1}{4}$. On $\Rea u=\frac{3}{4}+\eta$, for every $\eps>0$,
\begin{equation}\label{eq:low}
 T_\beta(u;M,D)\ll_{U,V,\eta,\eps}
 M^{\frac{7}{4}+\eps}D^{\frac{7}{2}+\eps}(1+|\Ima u|)^B.
\end{equation}
On $\Rea u=\frac{3}{2}+\eta$, for every $\eps>0$,
\begin{equation}\label{eq:high}
 T_\beta(u;M,D)\ll_{U,V,\eta,\eps}
 M^{\frac{1}{2}+\eps}D^{\frac{3}{2}+\eps}(1+|\Ima u|)^B.
\end{equation}
Consequently, on each fixed strip contained in $\frac{3}{4}<\sigma=\Rea u<\frac{3}{2}$, for every $\eps>0$,
\begin{equation}\label{eq:interpolation}
 T_\beta(u;M,D)\ll
 M^{3-\frac{5\sigma}{3}+\eps}D^{\frac{11}{2}-\frac{8\sigma}{3}+\eps}(1+|\Ima u|)^B.
\end{equation}
These estimates hold uniformly for $M,D\geq1$ and $\beta\in\{\beta_0,\beta_1\}$. In \eqref{eq:low} and \eqref{eq:high}, $B>0$ depends only on $\eta$ and $\eps$; in \eqref{eq:interpolation}, $B$ and the implied constant may also depend on the strip. The implied constants depend on $U,V$ only through their supports and upper bounds for finitely many of their derivatives.
\end{proposition}
\begin{proof}
On the line $\Rea u=\frac{3}{4}+\eta$, apply Proposition~\ref{prop:coprime} to the finite sum and take absolute values. The numbers of $m$ and $d$ are $O(M)$ and $O(D)$; the external factor $d$ contributes $D$, and the twist is $md^2$. This gives
\begin{equation}\label{eq:absolute-weighted-bound}
 T_\beta(u;M,D)\ll
 MD^2(MD^2)^{\frac{3}{2}-\Rea u+\eps}(1+|u|)^B,
\end{equation}
which implies \eqref{eq:low}.

On the line $\Rea u=\frac{3}{2}+\eta$, since $\Rea u>\frac{3}{2}$ and
\begin{equation}\label{eq:gauss-size}
 |g_K(\chi_c)|\leq\N(c)^{\frac{1}{2}},
\end{equation}
the sum over $c$ converges absolutely. Here and below, $\chi_c$ in a Dirichlet $L$-function denotes the restriction of the quartic residue character to rational integers, extended by zero at primes dividing $\N(c)$. Opening \eqref{eq:H}, we obtain
\begin{align}
 T_\beta(u;M,D)
 &=(u-\frac54)\zk(4u-3)
 \sum_{c\equiv\beta\ (4)}\frac{g_K(\chi_c)}{\N(c)^u}\notag\\
 &\quad\times\sum_{\substack{m,d\geq1\\(d,2m)=1}}
 \mu(d)d\,\bar\chi_c(m)\bar\chi_c(d)^2
 W_u\left(\frac{m}{M}\right)V\left(\frac{d}{D}\right).
 \label{eq:opened-weighted-sum}
\end{align}
The zero extension includes the coprimality condition in \eqref{eq:H}. Apply \eqref{eq:euler} with $\chi=\bar\chi_c$. Mellin inversion gives the explicit formula
\begin{align}
 T_\beta(u;M,D)
 &=(u-\frac54)\zk(4u-3)
 \sum_{c\equiv\beta\ (4)}\frac{g_K(\chi_c)}{\N(c)^u}
 \frac1{(2\pi i)^2}\int_{(2)}\int_{(2)}
 \widehat W_u(z)\widehat V(t+1)M^zD^{t+1}\notag\\
 &\quad\times\frac{L(z,\bar\chi_c)}{L^{(2)}(t,\bar\chi_c^2)}
 \prod_{p\ne2}\left(1+
 \frac{\bar\chi_c(p)^3p^{-z-t}}{1-\bar\chi_c(p)^2p^{-t}}\right)
 \dd t\,\dd z.
 \label{eq:double-mellin}
\end{align}
The factor $D^{t+1}$ occurs because of the weight $d$ in \eqref{eq:opened-weighted-sum}. We move both contours in \eqref{eq:double-mellin} to $\Rea z=\Rea t=\frac{1}{2}+\eps$. The only numerator pole crossed is at $z=1$ for a principal character; it is canceled by $1-z$ in \eqref{eq:weight}. The function $\frac{1}{L^{(2)}(t,\chi^2)}$ has no poles in this region under the hypothesis. For a principal character, this function has a zero at $t=1$.

The product in \eqref{eq:double-mellin} is normally convergent on the new contours. The bounds for $L(z,\chi)$ and $\frac{1}{L^{(2)}(t,\chi^2)}$, together with the finite Euler factors, contribute only a small power of the modulus and a polynomial in the vertical parameters. Choose that small power below the fixed margin $\Rea u-\frac{3}{2}$. The modulus sum remains absolutely convergent. Thus the shifted integral satisfies
\begin{equation}\label{eq:shifted-mellin-bound}
 T_\beta(u;M,D)\ll
 M^{\frac{1}{2}+\eps}D^{\frac{3}{2}+\eps}(1+|u|)^B
 \sum_{c\equiv\beta\ (4)}\N(c)^{\frac{1}{2}-\Rea u+\frac{\eta}{2}},
\end{equation}
where the conductor loss is chosen to be $\frac{\eta}{2}$. On $\Rea u=\frac{3}{2}+\eta$, the last sum has exponent $-1-\frac{\eta}{2}$ and converges. This proves \eqref{eq:high}. The factors $(u-\frac{5}{4})\zk(4u-3)$ contribute only polynomial vertical growth on this line.

Finally, apply the Phragm\'en--Lindel\"of principle to the holomorphic function $T_\beta(u;M,D)$ between the lines $\Rea u=\frac{3}{4}+\eta$ and $\Rea u=\frac{3}{2}+\eta$. Polynomial regularizing factors in $u$ take care of the vertical growth. As $\eta\to0^+$, the interpolated exponents of $M$ and $D$ tend to $3-\frac{5\sigma}{3}$ and $\frac{11}{2}-\frac{8\sigma}{3}$, respectively. On each fixed strip, choose $\eta$ and the losses in the estimates on these two lines sufficiently small in terms of $\eps$. This gives \eqref{eq:interpolation}. The Mellin estimates use only finitely many derivatives of $U,V$, which gives the asserted dependence of the implied constants.
\end{proof}

\subsubsection{Mellin inversion and continuation}
Fix $\delta\in\{0,1\}$ and set $\beta=\beta_\delta$ as in \eqref{eq:beta-parity}. We recover the unweighted series by a continuous Mellin identity. In a common region of absolute convergence, Fubini's theorem and the substitutions $x=\frac{m}{M}$ and $y=\frac{d}{D}$ yield
\begin{align}
 &\int_0^\infty\int_0^\infty
 T_\beta(u;M,D)M^{-v}D^{-2u}\frac{\dd M}{M}\frac{\dd D}{D}\notag\\
 &\quad=\widehat W_u(v)\widehat V(2u)
 (u-\frac54)\zk(4u-3)
 \bigl(C_{4,\delta}(u,v)+\mathbf1_{\delta=0}\zeta(v)\bigr).
 \label{eq:reconstruction}
\end{align}
The integrand vanishes when either $M$ or $D$ is below a fixed positive constant. Multiplying \eqref{eq:interpolation} by $M^{-v}D^{-2u}$ gives
\begin{equation}\label{eq:mellin-integrand-bound}
 |T_\beta(u;M,D)M^{-v}D^{-2u}|
 \ll M^{3-\frac{5\sigma}{3}-\Rea v+\eps}
 D^{\frac{11}{2}-\frac{14\sigma}{3}+\eps}(1+|u|)^B.
\end{equation}
The two exponents in \eqref{eq:mellin-integrand-bound} are negative for a sufficiently small $\eps$ precisely in the region
\begin{equation}\label{eq:uv-region}
 R_2=\{(u,v):\Rea(v+\frac53u)>3,\quad
 \frac{33}{28}<\Rea u<\frac32\}.
\end{equation}
The integrals in \eqref{eq:reconstruction} therefore converge normally in \eqref{eq:uv-region}.

At any prescribed point one can choose $U,V$ so that $\widehat U(v)$ and $\widehat V(2u)$ are nonzero. Dividing \eqref{eq:reconstruction} locally, and patching on overlaps by the identity theorem, continues the same series. The factors $1-v$ and $2-u-v$ are kept explicit. We use the change of variables
\begin{equation}\label{eq:change-variables}
 u=s+w,\qquad v=1-w,\qquad 1-v=w,\qquad 2-u-v=1-s.
\end{equation} The term $\mathbf1_{\delta=0}\zeta(v)$ is subtracted explicitly.

The estimate on the right hand line also gives a holomorphic continuation, with the displayed factors retained, in
\begin{equation}\label{eq:right-region}
 R_1=\{(u,v):\Rea u>\frac32,\ \Rea v>\frac12\}.
\end{equation} It meets the original absolute domain $\Rea u>\frac{3}{2}$, $\Rea v>1$. Interpolation between the lines $\Rea u=\frac{3}{4}+\eta$ and $\Rea u=\frac{3}{2}+\eta$, with $\eta>0$ sufficiently small, attaches \eqref{eq:uv-region} to this domain through a nonempty open overlap. The two representations therefore agree on a nonempty open set.

By \eqref{eq:change-variables}, the image of \eqref{eq:uv-region} in the original variables is
\begin{equation}\label{eq:new-region}
 S_3=\{(s,w):5\Rea s+2\Rea w>6,\quad
 \frac{33}{28}<\Rea(s+w)<\frac32\}.
\end{equation}
\subsubsection{The convex hull and vertical growth}
The boundary vertex of $S_3$ is $C=(\frac{17}{14},-\frac{1}{28})$. The region $S_1$ has limiting corner $A=(\frac{1}{2},1)$. The line through these points is
\begin{equation}\label{eq:hull-line}
 58\Rea s+40\Rea w=69.
\end{equation}
Apply Theorem~\ref{thm:tube} to the common holomorphic function obtained by multiplying by $s-1$, $s+w-\frac{5}{4}$ and the factors removing the finitely many gamma poles in the bounded real strips in question. Keep the same completion $\zk(4s+4w-3)$ during this step. Strict displacements of $A$ and $C$ lie in the source regions, so interior convex combinations are available. Only the portion of \eqref{eq:new-region} near $C$ is needed; all estimates are used within their stated strips.

In a neighborhood of $w=\frac{1}{2}$, the continued region contains the part needed below of
\begin{equation}\label{eq:final-region}
 S_4=\{(s,w):\Rea s>\frac12,\quad
 \Rea(s+w)>\frac{33}{28},\quad
 58\Rea s+40\Rea w>69\}.
\end{equation}
Throughout the needed hull,
\begin{equation}\label{eq:completion-halfplane}
 \Rea(4s+4w-3)>\frac{12}{7}>1.
\end{equation}
By \eqref{eq:completion-halfplane}, the completion can be divided out without introducing poles. The possible pole on $s+w=\frac{5}{4}$ remains explicit; on the central contour considered here it is strictly to the left.

\begin{lemma}\label{lem:growth}
The continued series has polynomial vertical growth on compact real subsets of the above domain away from its poles. In particular, for every fixed $\eps>0$, there is $B>0$ such that
\begin{equation}\label{eq:central-growth}
 |(s-1)A_4(s,\frac12)|\ll_\eps(1+|s|)^B,
 \qquad \frac{49}{58}+\eps\leq\Rea s\leq1+\eps.
\end{equation}
\end{lemma}
\begin{proof}
For the Mellin construction, modulate fixed smooth compactly supported $U,V$ by powers $x^{-it}$. Their Mellin transforms at the prescribed imaginary coordinates become nonzero fixed real Mellin integrals. The suprema of the absolute values of their derivatives of any fixed order grow at most polynomially in the modulation parameters. Propositions~\ref{prop:coprime} and \ref{prop:joint-bound}, followed by interpolation, therefore give polynomial bounds for the locally divided expression in \eqref{eq:reconstruction}. This avoids division by a globally small Mellin transform.

The initial Hecke representation has the same growth under the hypothesis. The gamma ratios in \eqref{eq:FE} have polynomial growth by Stirling's formula. To make the growth step explicit, restrict to the bounded real source strips used in the convex combination. Choose a real $C_0$ so that $s+C_0$ and $w+C_0$ have no zeros there, and choose $N$ large enough that
\begin{equation}\label{eq:majorant}
 h(s,w)=C(s+C_0)^N(w+C_0)^N
\end{equation}
dominates the absolute value of the common holomorphic function on the source union. Theorem~\ref{thm:tube-bound} preserves this domination on the convex hull. Finally, $\frac{1}{\zk(z)}$ is bounded in each fixed half plane $\Rea z\geq1+\eta$, where $\eta>0$. These statements prove the assertion.
\end{proof}

\subsection{Residue of \texorpdfstring{$A_4(s,w)$}{A4(s,w)} at \texorpdfstring{$s=1$}{s=1}}
For the density of the conductors, put
\begin{equation}\label{eq:density}
 \mathfrak c_4=\frac{\pi}{8}
 \prod_{p\equiv1\ (4)}\left(1-\frac{3}{p^2}+\frac{2}{p^3}\right)
 \prod_{p\equiv3\ (4)}\left(1-\frac{1}{p^2}\right)
\end{equation}
and
\begin{equation}\label{eq:constant}
 C_4=\mathfrak c_4\zeta(2)
 \prod_{p\equiv1\ (4)}\left(1-\frac{2}{p^2(p+2)}\right).
\end{equation}
These products converge absolutely and are positive.

We now compute the residue. In the initial expansion, the pole at $s=1$ comes from the terms with $m$ a fourth power. Indeed, the Kummer character $\chi_4^{(m)}$ is principal precisely when $m\in K^{\times4}$. A positive rational number which is a fourth power in $\Q(i)$ is a rational fourth power: a fourth root with equal nonzero real and imaginary parts has a negative fourth power, whereas roots on the coordinate axes give rational fourth powers. For positive integral $m$, the rational fourth root is integral. For such $m=a^4$, the conductor sum, with $n=1$ included, is
\begin{equation}\label{eq:fourth-power-conductors}
 1+\sideset{}{'}{\sum}_{\substack{n\ {\rm primary}\\(n,a)=1}}
 \frac1{\N(n)^s}
 =\prod_{\substack{p\equiv1\ (4)\\p\nmid a}}(1+2p^{-s}).
\end{equation}
The omitted term $q=1$ is entire and does not affect the residue. Factoring out $\zeta(s)L(s,\chi_{-4})$ gives
\begin{align}
 \prod_{p\equiv1\ (4)}(1+2p^{-s})
 &=\zeta(s)L(s,\chi_{-4})(1-2^{-s})\notag\\
 &\quad\times\prod_{p\equiv1\ (4)}(1-3p^{-2s}+2p^{-3s})
       \prod_{p\equiv3\ (4)}(1-p^{-2s}).
 \label{eq:conductor-euler}
\end{align}
The residue of \eqref{eq:conductor-euler} is $\mathfrak c_4$ in \eqref{eq:density}, since $L(1,\chi_{-4})=\frac{\pi}{4}$. Removing the primes dividing $a$ gives
\begin{equation}\label{eq:removed-prime-residue}
 \Res_{s=1}\prod_{\substack{p\equiv1\ (4)\\p\nmid a}}(1+2p^{-s})
 =\mathfrak c_4\prod_{\substack{p\mid a\\p\equiv1\ (4)}}\frac{p}{p+2}.
\end{equation}
Summing \eqref{eq:removed-prime-residue} over $a$, initially for $\Rea w>1$ and then by continuation, we obtain
\begin{align}
 \Res_{s=1}A_4(s,w)
 &=\mathfrak c_4\sum_{a\geq1}\frac1{a^{4w}}
              \prod_{\substack{p\mid a\\p\equiv1\ (4)}}\frac{p}{p+2}\notag\\
 &=\mathfrak c_4\zeta(4w)
        \prod_{p\equiv1\ (4)}\left(1-\frac{2p^{-4w}}{p+2}\right).
 \label{eq:residue}
\end{align}
This is the quartic specialization of the residue calculation in [Section 3.4 of~\citen{GZ}], with the squarefree Euler denominator interpreted as $L_K(2s,(\chi_4^{(m)})^2)$. At $w=\frac{1}{2}$, \eqref{eq:residue} is exactly $C_4$ in \eqref{eq:constant}.

\subsection{Completion of the proof}
\begin{proof}[Proof of the quartic case of Theorem~\ref{thm:main}]
We apply \eqref{eq:moment-integral} with $w=\frac{1}{2}$. For $c>1$, the stated hypothesis gives
\begin{equation}\label{eq:central-integral}
 \sum_{\substack{q>1\\(q,2)=1}}\ \sumstar_{\substack{\chi\bmod q\\\chi^4=\chi_0}}
 L(\frac12,\chi)\Phi\left(\frac{q}{Q}\right)
 =\frac1{2\pi i}\int_{(c)}A_4(s,\frac12)Q^s\widehat\Phi(s)\dd s.
\end{equation}
At $w=\frac{1}{2}$, the boundary \eqref{eq:hull-line} is $\Rea s=\frac{49}{58}$. We move the contour in \eqref{eq:central-integral} to $\Rea s=\frac{49}{58}+\eps$. The only pole crossed is at $s=1$, whose residue is given by \eqref{eq:residue}. Hence
\begin{align}
 &\sum_{\substack{q>1\\(q,2)=1}}\ \sumstar_{\substack{\chi\bmod q\\\chi^4=\chi_0}}
 L(\frac12,\chi)\Phi\left(\frac{q}{Q}\right)\notag\\
 &\qquad=C_4\widehat\Phi(1)Q+
 \frac1{2\pi i}\int_{(\frac{49}{58}+\eps)}
 A_4(s,\frac12)Q^s\widehat\Phi(s)\dd s.
 \label{eq:shifted-moment}
\end{align}
The possible second pole at $s=\frac{3}{4}$ lies to the left of the new contour. Integration by parts in \eqref{eq:mellin-transform} gives, for every $N\geq0$ and on each fixed real strip,
\begin{equation}\label{eq:mellin-decay}
 \widehat\Phi(s)\ll_{\Phi,N}(1+|\Ima s|)^{-N}.
\end{equation}
By \eqref{eq:central-growth} and \eqref{eq:mellin-decay},
\begin{align}
 &\left|\frac1{2\pi i}\int_{(\frac{49}{58}+\eps)}
 A_4(s,\frac12)Q^s\widehat\Phi(s)\dd s\right|\notag\\
 &\qquad\ll_{\Phi,\eps}Q^{\frac{49}{58}+\eps}
 \int_{\mathbb R}(1+|t|)^B
 |\widehat\Phi(\frac{49}{58}+\eps+it)|\dd t
 \ll_{\Phi,\eps}Q^{\frac{49}{58}+\eps}.
 \label{eq:contour-error}
\end{align}
The same bounds justify the vanishing of the horizontal integrals in the contour shift. Combining \eqref{eq:shifted-moment} and \eqref{eq:contour-error} proves the case $j=4$ of \eqref{eq:main}.
\end{proof}

\section{The sextic first moment}\label{sec:sextic}
\subsection{Sextic characters and Gauss sums}
In this section let $K_6=\Q(\sqrt{-3})$ and $\OO_{K_6}=\Z[\omega]$, where $\omega=(-1+\sqrt{-3})/2$. We use the E-primary generators of [Section 2.1 of~\citen{GZ}]. Explicitly, $n=a+b\omega$ is E-primary if $(n,6)=1$, $n\equiv\pm1\pmod3$, and
\begin{equation}\label{eq:E-primary}
 \begin{cases}
 a+b\equiv1\pmod4,&2\mid b,\\
 b\equiv1\pmod4,&2\mid a,\\
 a\equiv3\pmod4,&2\nmid ab.
 \end{cases}
\end{equation}
Every ideal coprime to $6$ has a unique such generator. Products of E-primary elements are E-primary. The sextic residue symbol is extended by zero when its arguments are not coprime, and we write
\begin{equation}\label{eq:sextic-symbols}
 \chi_{6,n}(m)=\left(\frac{m}{n}\right)_6,
 \qquad \chi_6^{(m)}(c)=\left(\frac{m}{c}\right)_6.
\end{equation}
As in the quartic case, the first character is restricted to rational integers when it occurs in a Dirichlet $L$-function. The second defines a finite order Hecke character on ideals coprime to $6m$ by evaluation on their unique E-primary generators. This is multiplicative because products of E-primary generators are E-primary. Its modulus includes the fixed conditions at the primes above $2$ and $3$, and it is extended by zero at the original prime divisors of $m$. All factors at the primes above $2$ and $3$ are retained or explicitly omitted as indicated below.

By [Lemma 2.2 of~\citen{GZ}], the sextic family in Theorem~\ref{thm:main} is parameterized by squarefree E-primary $n\ne1$ without a rational prime divisor, with $q=\N(n)$. We denote this restriction by $\sum'$. Define
\begin{equation}\label{eq:sextic-A}
 A_6(s,w)=\sum_{\substack{q>1\\(q,6)=1}}
 \sumstar_{\substack{\chi\bmod q\\\chi^6=\chi_0}}\frac{L(w,\chi)}{q^s}
 =\sideset{}{'}{\sum}_{n\ \mathrm{E\text{-}primary}}
 \frac{L(w,\chi_{6,n})}{\N(n)^s}.
\end{equation}
The sum initially converges absolutely for $\Rea s,\Rea w>1$.

We use the additive character and Gauss sums
\begin{equation}\label{eq:sextic-gauss}
 e_{K_6}(z)=\exp\left(2\pi i\frac{z-\bar z}{\sqrt{-3}}\right),\qquad
 g_{K_6}(r,\chi_{6,n})=\sum_{x\bmod n}\left(\frac{x}{n}\right)_6
 e_{K_6}\left(\frac{rx}{n}\right).
\end{equation}
We abbreviate $g_{K_6}(1,\chi_{6,n})$ to $g_{K_6}(\chi_{6,n})$. By [Lemma 2.5 of~\citen{GZ}], its absolute value is $\sqrt{\N(n)}$ for squarefree $n$, and it is zero otherwise. We will need its value, not just its absolute value, when the denominator is rational.

\begin{lemma}\label{lem:sextic-rational-gauss}
Let $d$ be a positive squarefree rational integer with $(d,6)=1$, represented by its E-primary associate. Then
\begin{equation}\label{eq:sextic-rational-gauss}
 g_{K_6}(\chi_{6,d})=d.
\end{equation}
For coprime rational integers $m,d$ with $(d,6)=1$, the sextic symbol $(m/d)_6$ equals $1$.
\end{lemma}
\begin{proof}
First let $p\nmid6$ be a rational prime and put $\vartheta(x)=(x/p)_6$ on $\OO_{K_6}/(p)$. This is a nonprincipal multiplicative character, extended by zero off the units, whose restriction to $\mathbb F_p^*$ is trivial. For split $p$, the restrictions from the two conjugate factors cancel. For inert $p$, $6\mid p+1$, which gives the same restriction. These statements also prove the assertion about rational symbols, by multiplicativity.

Write $x=a+b\omega$ and let $S_b=\sum_{a\bmod p}\vartheta(a+b\omega)$. Then $S_0=p-1$. Scaling by a nonzero rational $b$ gives $S_b=S_1$ for $b\ne0$. The complete multiplicative character sum is zero, so $S_1=-1$. Since the additive character at rational denominator $p$ is $e(b/p)$, we obtain
\begin{equation}\label{eq:sextic-rational-prime}
 g_{K_6}(\chi_{6,p})=p-1-\sum_{b\ne0}e(b/p)=p.
\end{equation}
This argument applies to both the split algebra and the inert field. The Chinese Remainder Theorem gives multiplicativity in squarefree rational denominators, since the cross symbols are trivial. Passing to the E-primary rational associate changes only a sign, on which the rational restriction is trivial. This proves \eqref{eq:sextic-rational-gauss}.
\end{proof}

\subsection{The functional equation and the rational divisor sum}
Put $c_6=-9$. For $a\in\{1,-1\}$, define
\begin{align}
 C_{6,a}(u,v)
 &=\sum_{m\geq1}\frac1{m^v}
 \sideset{}{'}{\sum}_{n\ \mathrm{E\text{-}primary}}
 \frac{g_{K_6}(a c_6m,\chi_{6,n})}{\N(n)^u},\label{eq:sextic-C}\\
 H_{6,a}(u,r)
 &=\sum_{\substack{c\ \mathrm{E\text{-}primary}\\(c,6r)=1}}
 \frac{g_{K_6}(a c_6r,\chi_{6,c})}{\N(c)^u}.
 \label{eq:sextic-H}
\end{align}
The sign $a$ has the same role as in $C_{j,a}$ in [Section 3.3 of~\citen{GZ}]. It does not denote parity. The factor $c_6$ is included in the frequency in \eqref{eq:sextic-H}; this accounts for the difference from the argument of $H_{j,a}$ in that paper.

By [Lemma 2.6 of~\citen{GZ}],
\begin{equation}\label{eq:sextic-tau}
 \tau(\chi_{6,n})=
 \overline{\left(\frac{c_6}{n}\right)_6}g_{K_6}(\chi_{6,n}).
\end{equation}
Thus $g_{K_6}(a c_6m,\chi_{6,n})=\chi_{6,n}(a)\overline{\chi_{6,n}(m)}\tau(\chi_{6,n})$. With $\gamma_\delta$ as in \eqref{eq:dirichlet-FE}, the ordinary Dirichlet functional equation gives the exact identity
\begin{align}
 A_6(s,w)
 &=\frac{\gamma_0(w)}2
 \bigl(C_{6,1}(s+w,1-w)+C_{6,-1}(s+w,1-w)\bigr)\notag\\
 &\quad-\frac{i\gamma_1(w)}2
 \bigl(C_{6,1}(s+w,1-w)-C_{6,-1}(s+w,1-w)\bigr).
 \label{eq:sextic-FE}
\end{align}
In particular, the odd parity factor $i^{-1}$ is retained here. The convention for $\gamma_\delta$ has not changed from the quartic calculation.

M\"obius inversion and [Lemma 2.5 of~\citen{GZ}] give, in an absolute domain,
\begin{equation}\label{eq:sextic-mobius}
 C_{6,a}(u,v)+\zeta(v)=
 \sum_{m\geq1}\frac1{m^v}
 \sum_{\substack{d\geq1\\(d,6m)=1}}
 \frac{\mu(d)d}{d^{2u}}H_{6,a}(u,md^4).
\end{equation}
Indeed, the Gauss sum identity introduces $d^{6-2}$ in the frequency, and Lemma~\ref{lem:sextic-rational-gauss} evaluates the outside coefficient as $d$. The rational symbol is $1$ when $(m,d)=1$ and zero otherwise. The added term $\zeta(v)$ is the contribution of $n=1$ for each sign $a$.

For some continuation arguments it is convenient to write
\begin{equation}\label{eq:sextic-unit}
 A_6^+(s,w)=A_6(s,w)+\zeta(w),\qquad
 C_{6,a}^+(u,v)=C_{6,a}(u,v)+\zeta(v).
\end{equation}
The added terms satisfy \eqref{eq:sextic-FE} by the functional equation of $\zeta$. Their contribution to the smoothed moment is $\zeta(\frac12)\Phi(1/Q)$, which vanishes for all sufficiently large $Q$.

\subsection{An estimate for the sextic Gauss series}
The next input is already available for sextic characters. We use it only in $\Rea u>1$.

\begin{lemma}\label{lem:sextic-gauss-bound}
The function $(u-\frac76)H_{6,a}(u,r)$ is holomorphic in $\Rea u>1$. On fixed strips contained in $1<\sigma=\Rea u<\frac32$, for every $\eps>0$,
\begin{equation}\label{eq:sextic-gauss-bound}
 (u-\tfrac76)H_{6,a}(u,r)
 \ll_{\eps}r^{\frac32-\sigma+\eps}(1+|u|)^B.
\end{equation}
The estimate holds uniformly for positive integers $r$ and the two signs $a$. Versions on the lines $\sigma=1+\eta$ and $\sigma=\frac32+\eta$ hold with arbitrarily small losses, for sufficiently small fixed $\eta>0$.
\end{lemma}
\begin{proof}
Apply [Lemma 2.9 of~\citen{GZ}] to the coprime Gauss series with frequency $a c_6r$. The norm of this frequency is $81r^2$, and the exponent in that lemma is $(\frac32-\sigma)/2$, up to an arbitrarily small loss. This proves \eqref{eq:sextic-gauss-bound}. The polynomial in the cited lemma contains the factors $u-\frac76$ and $u-\frac56$. The latter is bounded away from zero in the strips used here and can be divided out. The bounded conductor characters required to impose the E-primary conditions do not alter uniformity in $r$.
\end{proof}

In particular,
\begin{equation}\label{eq:sextic-twist-bound}
 (u-\tfrac76)H_{6,a}(u,md^4)
 \ll m^{\frac32-\sigma+\eps}d^{6-4\sigma+\eps}(1+|u|)^B.
\end{equation}
This dependence on $d$ is different from that in the quartic case and will give different interpolation exponents.

\subsection{The joint estimate}
For a Dirichlet character $\chi$, extended by zero at every prime of its stated modulus, write $L^{(6)}(t,\chi)$ for the Euler product with the factors at $2$ and $3$ omitted. The analogue of Lemma~\ref{lem:euler} is
\begin{equation}\label{eq:sextic-euler}
 \sum_{\substack{m,d\geq1\\(m,d)=1,\ (d,6)=1}}
 \frac{\mu(d)\chi(m)\chi(d)^4}{m^zd^t}
 =\frac{L(z,\chi)}{L^{(6)}(t,\chi^4)}
 \prod_{p\nmid6}\left(1+
 \frac{\chi(p)^5p^{-z-t}}{1-\chi(p)^4p^{-t}}\right).
\end{equation}
Initially $\Rea z,\Rea t>1$. The correction product converges normally for $\Rea t>0$ and $\Rea(z+t)>1$. To verify the identity, the local sum at $p\nmid6$ is
\begin{equation}\label{eq:sextic-euler-local}
 \frac1{1-\chi(p)p^{-z}}-\chi(p)^4p^{-t}.
\end{equation}
Factoring out $(1-\chi(p)^4p^{-t})/(1-\chi(p)p^{-z})$ gives \eqref{eq:sextic-euler}. At $2$ and $3$, only the $m$ sum remains. At any prime of the original modulus the character is zero and the local factor is one. We keep this zero extension even if the primitive conductor is smaller. In the application, $\chi=\bar\chi_{6,c}$ restricted to rational integers, so $\chi^4$ has order dividing three.

Use the weights $U,V,W_u$ in \eqref{eq:weight}, and define
\begin{equation}\label{eq:sextic-weighted-sum}
 T_{6,a}(u;M,D)=(u-\tfrac76)
 \sum_{\substack{m,d\geq1\\(d,6m)=1}}
 \mu(d)d H_{6,a}(u,md^4)
 W_u\left(\frac mM\right)V\left(\frac dD\right).
\end{equation}
By Lemma~\ref{lem:sextic-gauss-bound}, this finite sum is holomorphic for $\Rea u>1$.

\begin{proposition}\label{prop:sextic-joint}
For arbitrarily small positive displacements of the indicated lines, we have
\begin{align}
 T_{6,a}(u;M,D)&\ll M^{\frac32+\eps}D^{4+\eps}(1+|u|)^B,
 &&\Rea u=1+\eta,\label{eq:sextic-low}\\
 T_{6,a}(u;M,D)&\ll M^{\frac12+\eps}D^{\frac32+\eps}(1+|u|)^B,
 &&\Rea u=\tfrac32+\eta.\label{eq:sextic-high}
\end{align}
Consequently, on each fixed strip contained in $1<\sigma=\Rea u<\frac32$,
\begin{equation}\label{eq:sextic-interpolation}
 T_{6,a}(u;M,D)\ll
 M^{\frac72-2\sigma+\eps}D^{9-5\sigma+\eps}(1+|u|)^B.
\end{equation}
The constants are uniform for $M,D\geq1$ and depend on the weights only through their supports and finitely many smooth seminorms.
\end{proposition}
\begin{proof}
For \eqref{eq:sextic-low}, apply \eqref{eq:sextic-twist-bound} termwise. There are $O(MD)$ terms, and the outside factor $d$ contributes $D$. At $\sigma=1$ the resulting powers are $M^{3/2}D^4$.

For \eqref{eq:sextic-high}, open $H_{6,a}$ in \eqref{eq:sextic-weighted-sum}. With $\xi_c=\bar\chi_{6,c}$ on rational integers, including its original zero extension, the inner sum is
\begin{equation}\label{eq:sextic-inner-sum}
 \sum_{\substack{m,d\geq1\\(d,6m)=1}}
 \mu(d)d\,\xi_c(m)\xi_c(d)^4W_u(m/M)V(d/D).
\end{equation}
The outside multiplier is $g_{K_6}(a c_6,\chi_{6,c})/\N(c)^u$ and is independent of $m,d$. By \eqref{eq:sextic-euler}, Mellin inversion expresses \eqref{eq:sextic-inner-sum} as
\begin{align}
 \frac1{(2\pi i)^2}\int_{(2)}\int_{(2)}
 &\widehat W_u(z)\widehat V(t+1)M^zD^{t+1}
 \frac{L(z,\xi_c)}{L^{(6)}(t,\xi_c^4)}\notag\\
 &\times\prod_{p\nmid6}\left(1+
 \frac{\xi_c(p)^5p^{-z-t}}{1-\xi_c(p)^4p^{-t}}\right)\dd t\,\dd z.
 \label{eq:sextic-double-mellin}
\end{align}
Move both contours to $\Rea z=\Rea t=\frac12+\eps$. The numerator pole at $z=1$ for a principal character is canceled by $1-z$ in \eqref{eq:weight}. Under GRH the reciprocal $L$-function has no poles in the new half plane, and the correction product is normally convergent. The factors at the original modulus cost at most $\N(c)^\eps$, including when the primitive conductor drops. Since $|g_{K_6}(a c_6,\chi_{6,c})|\leq\sqrt{\N(c)}$, the sum over $c$ still converges absolutely when the conductor loss is smaller than the fixed margin $\Rea u-\frac32$. The scale factors in \eqref{eq:sextic-double-mellin} therefore give \eqref{eq:sextic-high}.

The Phragm\'en--Lindel\"of principle applied to the same holomorphic finite sum between these lines gives \eqref{eq:sextic-interpolation}. The limiting affine exponents are $\frac72-2\sigma$ and $9-5\sigma$. As in Proposition~\ref{prop:joint-bound}, the small line displacements are absorbed into $\eps$ on each fixed strip. Polynomial regularizing factors control vertical growth.
\end{proof}

\subsection{Mellin inversion and the regions of continuation}
The analogue of \eqref{eq:reconstruction} is the exact identity
\begin{align}
 &\int_0^\infty\int_0^\infty
 T_{6,a}(u;M,D)M^{-v}D^{-2u}\frac{\dd M}{M}\frac{\dd D}{D}\notag\\
 &\qquad=\widehat W_u(v)\widehat V(2u)
 (u-\tfrac76)\bigl(C_{6,a}(u,v)+\zeta(v)\bigr).
 \label{eq:sextic-reconstruction}
\end{align}
It follows from Fubini's theorem in an absolute domain. The small scale portions vanish below fixed positive constants. By \eqref{eq:sextic-interpolation}, the large scale integrals converge normally when
\begin{equation}\label{eq:sextic-source}
 \Rea u>\frac97,\qquad \Rea(v+2u)>\frac72,
 \qquad 1<\Rea u<\frac32.
\end{equation}
Indeed, the powers of $M$ and $D$ in the integrand are $\frac72-2\Rea u-\Rea v+\eps$ and $9-7\Rea u+\eps$.

At any prescribed point, choose $U,V$ with nonzero Mellin transforms there. Local division in \eqref{eq:sextic-reconstruction}, retaining the factors $(1-v)(2-u-v)$, and the identity theorem give continuation of the same series. The high estimate also gives the open region $\Rea u>\frac32$, $\Rea v>\frac12$, which meets the original absolute domain $\Rea u>\frac32$, $\Rea v>1$. The displaced interpolation strips meet this region. Thus all the representations agree on open overlaps. Modulating the weights as in the proof of Lemma~\ref{lem:growth} gives polynomial vertical bounds for the locally divided expressions.

For completeness, the initial expansion is
\begin{align}
 A_6(s,w)+\zeta(w)
 &=\sum_{m\geq1}\frac1{m^w}
 \frac{L_{K_6}^{(6)}(s,\chi_6^{(m)})}
 {L_{K_6}^{(6)}(2s,(\chi_6^{(m)})^2)}\notag\\
 &\quad\times\sum_{\substack{d\geq1\\(d,6m)=1}}\frac{\mu(d)}{d^{2s}}
 \prod_{\varpi\mid d}\left(1+
 \frac{\chi_6^{(m)}(\varpi)}{\N(\varpi)^s}\right)^{-1}.
 \label{eq:sextic-initial}
\end{align}
Here $\chi_6^{(m)}$ is the Hecke character defined using E-primary generators in \eqref{eq:sextic-symbols}, as in [Section 3.2 of~\citen{GZ}]. Its finite conditions at the primes above $6$ are part of its modulus. The factors above $6$ are omitted on both sides of each Euler product. The derivation and the estimates are the same as in \eqref{eq:first-mobius}--\eqref{eq:first-bound}, with $K_6$ in place of $K$. In particular, the denominator is evaluated at $2s$ and is absolutely convergent for $\Rea s>\frac12$. Hecke GRH gives the first region
\begin{equation}\label{eq:sextic-S1}
 S_{1,6}=\{(s,w):\Rea s>\tfrac12,\ \Rea w>1\},
\end{equation}
with the first main pole at $s=1$ retained. The conductor bound \eqref{eq:conductor-bound} gives the connecting region
\begin{equation}\label{eq:sextic-S2}
 S_{2,6}=\{(s,w):\Rea s>1,\ \Rea(s+w)>\tfrac32\},
\end{equation}
since there are at most $2^{\omega(q)}$ characters in the sextic family at conductor $q$. This region meets $S_{1,6}$ for large $s,w$, and meets the dual absolute domain for $\Rea w<0$ and sufficiently large $\Rea s$.

Under $u=s+w$, $v=1-w$, the new region \eqref{eq:sextic-source} becomes
\begin{equation}\label{eq:sextic-S3}
 S_{3,6}=\{(s,w):\Rea(s+w)>\tfrac97,
 \ \Rea(2s+w)>\tfrac52,\ \Rea(s+w)<\tfrac32\}.
\end{equation}
Its limiting corner is $(\frac{17}{14},\frac1{14})$, while the limiting corner of $S_{1,6}$ is $(\frac12,1)$. The line joining them is
\begin{equation}\label{eq:sextic-hull}
 26\Rea s+20\Rea w=33.
\end{equation}
Before applying Theorems~\ref{thm:tube} and \ref{thm:tube-bound}, multiply by the elementary factors removing the poles. For the enlarged series $A_6^+$, these include
\begin{equation}\label{eq:sextic-clearing}
 w(w-1)(s-1)(s+w-\tfrac76).
\end{equation}
The factor $w-1$ removes the pole of the added $\zeta(w)$ term. Any further gamma poles in the bounded real strips are removed as well. Strict displacements of the two corners give interior points of the convex hull. We obtain continuation to a neighborhood of every point on $w=\frac12$ with $\Rea s>\frac{23}{26}$, together with polynomial vertical bounds after clearing $s=1$. The factors $w$ and $w-1$ in \eqref{eq:sextic-clearing} are nonzero on this central line.

\subsection{The residue and the contour shift}
Let $\chi_{-3}$ be the primitive quadratic Dirichlet character modulo $3$. Define
\begin{align}
 \mathfrak c_6&=\frac{2\pi}{9\sqrt3}
 \prod_{p\equiv1\ (3)}\left(1-\frac3{p^2}+\frac2{p^3}\right)
 \prod_{p\equiv2\ (3)}\left(1-\frac1{p^2}\right),\label{eq:sextic-density}\\
 C_6&=\mathfrak c_6\zeta(3)
 \prod_{p\equiv1\ (3)}\left(1-\frac2{p^3(p+2)}\right).
 \label{eq:sextic-constant}
\end{align}
The second product in \eqref{eq:sextic-density} includes $p=2$. All products converge absolutely and are positive.

In the initial expansion, the pole at $s=1$ comes from positive rational $m$ which are sixth powers in $K_6$, and these are precisely rational sixth powers. To see this, a sixth root $z\in K_6$ of a positive rational number lies on a ray with angle an integer multiple of $\pi/3$. Multiplication by a unit makes it real, and $K_6\cap\mathbb R=\Q$. If $m$ is integral, its rational sixth root is integral. Thus we may write $m=b^6$.

For such $m$, the conductor Euler product, with the unit included, is
\begin{equation}\label{eq:sextic-conductor-product}
 \prod_{\substack{p\equiv1\ (3)\\p\nmid b}}(1+2p^{-s}).
\end{equation}
Without the restriction $p\nmid b$, this factors as
\begin{align}
 \prod_{p\equiv1\ (3)}(1+2p^{-s})
 &=\zeta(s)L(s,\chi_{-3})(1-3^{-s})\notag\\
 &\quad\times\prod_{p\equiv1\ (3)}(1-3p^{-2s}+2p^{-3s})
 \prod_{p\equiv2\ (3)}(1-p^{-2s}).
 \label{eq:sextic-conductor-factorization}
\end{align}
Since $L(1,\chi_{-3})=\pi/(3\sqrt3)$, its residue is \eqref{eq:sextic-density}. Omitting the primes dividing $b$ and summing over $b$ gives
\begin{align}
 \Res_{s=1}A_6(s,w)
 &=\mathfrak c_6\sum_{b\geq1}\frac1{b^{6w}}
 \prod_{\substack{p\mid b\\p\equiv1\ (3)}}\frac{p}{p+2}\notag\\
 &=\mathfrak c_6\zeta(6w)
 \prod_{p\equiv1\ (3)}\left(1-\frac{2p^{-6w}}{p+2}\right).
 \label{eq:sextic-residue}
\end{align}
This is the $j=6$ specialization of [Section 3.4 of~\citen{GZ}], written as a rational prime Euler product. At $w=\frac12$ it equals $C_6$ in \eqref{eq:sextic-constant}.

\begin{proof}[Proof of the sextic case of Theorem~\ref{thm:main}]
Mellin inversion gives, for $c>1$,
\begin{equation}\label{eq:sextic-moment-integral}
 \sum_{\substack{q>1\\(q,6)=1}}\sumstar_{\substack{\chi\bmod q\\\chi^6=\chi_0}}
 L(\tfrac12,\chi)\Phi(q/Q)
 =\frac1{2\pi i}\int_{(c)}A_6(s,\tfrac12)Q^s\widehat\Phi(s)\dd s.
\end{equation}
Move the contour to $\Rea s=\frac{23}{26}+\eps$. The only pole crossed is at $s=1$, and \eqref{eq:sextic-residue} gives its contribution as $C_6\widehat\Phi(1)Q$. The possible pole $s+w=\frac76$ meets the central line at $s=\frac23$ and lies strictly to the left of the new contour. Polynomial vertical growth and \eqref{eq:mellin-decay} show that the new integral is $O_{\Phi,\eps}(Q^{\frac{23}{26}+\eps})$ and that the horizontal integrals vanish. If $A_6^+$ is used in the continuation, subtract the unit term in \eqref{eq:sextic-unit}; it vanishes in the smoothed moment for sufficiently large $Q$. This proves the remaining case of \eqref{eq:main}.
\end{proof}

\end{document}